\documentclass[titlepage,11pt]{article}

\usepackage[pdftex]{graphicx}

\usepackage[myheadings]{fullpage}
\usepackage{pmetrika}
\usepackage{pmbib}
\usepackage{amsmath}
\usepackage{amsfonts}
\usepackage{amssymb}
\usepackage{mathtools}
\usepackage{cite}
\usepackage{apacite}
\usepackage{array}
\usepackage{enumitem}
\usepackage[nolists]{endfloat}

\usepackage{submit}

\newcommand\numberthis{\addtocounter{equation}{1}\tag{\theequation}}

\DeclareMathOperator\Var{Var}

\DeclareMathOperator\E{\mathbb{E}}

\DeclareMathOperator*{\dd}{\mathrm{d}}
\DeclareMathOperator*{\argmin}{arg\,min}

\begin{document}

\begin{titlepage}

\title{Optimal Designs for the Generalized Partial Credit Model}


\author{Paul-Christian B\"urkner$^1$, Rainer Schwabe$^2$, Heinz Holling$^1$}

\bigskip

\affil{$^1$ Institute of Psychology, University of M\"unster, Germany}
\affil{$^2$ Institute of Mathematical Stochastics, University of Magdeburg, Germany}

\vspace{\fill}\centerline{\today}\vspace{\fill}

\linespacing{1}

\end{titlepage}

\setcounter{page}{2}
\vspace*{2\baselineskip}

\RepeatTitle{Optimal Designs for the Generalized Partial Credit Model}\vskip3pt

\linespacing{1.5}
\abstracthead
\begin{abstract}
Analyzing ordinal data becomes increasingly important in psychology, especially in the context of item response theory. The generalized partial credit model (GPCM) is probably the most widely used ordinal model and finds application in many large scale educational assessment studies such as PISA. In the present paper, optimal test designs are investigated for estimating persons' abilities with the GPCM for calibrated tests when item parameters are known from previous studies. We will derive that local optimality may be achieved by assigning non-zero probability only to the first and last category independently of a person's ability. That is, when using such a design, the GPCM reduces to the dichotomous 2PL model. Since locally optimal designs require the true ability to be known, we consider alternative Bayesian design criteria using weight distributions over the ability parameter space. For symmetric weight distributions, we derive necessary conditions for the optimal one-point design of two response categories to be Bayes optimal. Furthermore, we discuss examples of common symmetric weight distributions and investigate, in which cases the necessary conditions are also sufficient. Since the 2PL model is a special case of the GPCM, all of these results hold for the 2PL model as well. 

\begin{keywords}
optimal design; Bayesian design; partial credit model; 2PL model; Rasch model; item response theory.
\end{keywords}
\end{abstract}\vspace{\fill}\pagebreak

\section{1. Introduction}

Item response theory (IRT) provides a flexible and powerful approach for designing and analyzing data of psychological tests. While IRT very often deals with dichotomous responses (e.g., categorized as either $1$ 'right' or $0$ 'wrong'), there are many situations where a more advanced scoring is to be advised. Consider, for instance, a complicated task in which multiple intermediate results have to be obtained to correctly solve the whole task. Scoring such a task in a dichotomous way comes with substantial loss of information as someone solving all intermediate questions failing only in the last step would receive the same score as someone who did not even manage to take the first step. A natural solution is to give persons \emph{partial credit} for the intermediate results. The obtained response is no longer dichotomous but ordinal ranging from $0$ (nothing correct), over $1$ (first step correct) to $J$ (all $J$ steps correct). An IRT model dealing with such ordinal responses is the Partial Credit Model (PCM; \citeNP{andersen1973, andrich1978a}) or its extention, the generalized partial credit model (GPCM;  \citeNP{muraki1992, muraki1993}). It is a mathematically convenient generalization of the dichotomous 2PL model (which is in turn a generalization of the Rasch model) for more than two possible response categories. As such, it combines the ability of an examinee with the item difficulty and discrimination and maps them to the probability of each response category. 

Although alternative ordinal models exist and are commonly applied in practice as well -- most notably the graded response model \cite{walker1967, vanderark2001} and the sequential model \cite{tutz1990, tutz2000} -- the GPCM is probably the most frequently used ordinal model in psychological research. Among others, it has been applied in many large scale educational assessment studies such as PISA, TIMSS, PIRLS, or NAEP \cite{oecd2017, martin2000, allen2001}. In addition, there is a great body of literature discussing application of the GPCM in large scale assessment studies from a more technical perspective (e.g., see \citeNP{mazzeo2014, vondavier2010, vondavier2013}). These studies were mostly concerned with comparability of large scale results as well as model fitting itself, while to our knowledge optimizing the statistical properties of the resulting parameter estimates has not yet been systematically investigated yet. We believe that given the GPCM's broad application in practice, it is of great relevance to investigate under which conditions the GPCM performs best. 

Applying principles of optimal design \cite{berger2009, atkinson2007, holling2013} to IRT models can lead to considerable efficiency gains and reduce costs of test administration by reducing the number of items and / or subjects required for achieving the same level of precision \cite{holling2016}. In IRT, one typically distinguishes between two types of optimal design problems. The first is about selecting items with specific properties for the efficient estimation of person parameters, while the second is about selecting persons with specific abilities for the efficient estimation of item parameters. These problems are referred to as optimal test design and optimal sampling design, respectively \cite{holling2016}, with the former being the more common one \cite{vanderlinden2006}. Optimal test designs are relevant in so called calibrated tests, in which item parameters were estimated in prior studies and can thus be considered as known with reasonable precision. The design question is to select the items from the calibrated item pool most appropriate for the participating persons. Optimal test designs are also applicable in tests making use of automatic item generation (e.g., see \citeNP{geerlings2011}). In this case, the item parameters can be inferred from the set of rules, which were used to create the items, provided that the influence of the rules on the item parameters were investigated before.

Most IRT models, including the GPCM, are non-linear in the parameters so that the related optimal designs are dependent on the parameters to be estimated. Such optimal designs that are optimal only for certain parameter values -- but not for others -- are called locally optimal designs. Among others, locally optimal designs for the GPCM will be investigated in the present paper. Despite being relevant for theoretical reasons and often being a prerequisite more advanced optimal designs, locally optimal designs are of limited practical use themselves as one has to guess the true value of the parameter before constructing the design. If the guess is far away from the truth, such locally optimal design might become relatively inefficient.

There are different possibilities to overcome the problem of locally optimal designs. One solution considered in the present paper are so called Bayesian optimal designs. Instead of trying to find an optimal design for a specific parameter value, one specifies a weight distribution over the set of possible parameters to express the uncertainty in the subject or item parameters. For instance, when generating a test design, one might assume that a subject's ability will be somewhere between $-2$ and $2$ with all abilities in between being equally likely. This would translate into a uniform distribution over the interval $[-2, 2]$. A design that is optimal given the specified uncertainty in the parameter values is then called Bayes optimal. 

The paper is structured as follows. In Section 2 and 3, the GPCM is introduced in detail and locally optimal designs are investigated. In Section 4, we introduce and discuss Bayesian optimal designs in relation to the GPCM. Examples for common weight distributions are given in Section 5. We end with a discussion of the findings in Section 6. All proofs are provided in the Appendix.

\section{2. The Model} 

In a test situation with ordinal items, the response of a person $p$ on item $i$ is symbolized by $Y_{pi}$, which can take on one of $J+1$ values ranging from $0$ to $J$. In the GPCM, the probability $\pi_{pij}$ that person $p$ achieves category $j \in \{0, ..., J\}$ on item $i$ is given by (cf. \citeNP{muraki1992, muraki1993}):

\begin{equation}
\label{PCM}
\pi_{pij} := \pi_j(\theta_p, \boldsymbol{\tau_i}, \boldsymbol{\alpha_i}) := P(Y_{pi} = j; \theta_p, \boldsymbol{\tau_i}, \boldsymbol{\alpha_i}) := 
 \frac{\exp \left(\sum_{s=1}^j \alpha_{is} (\theta_p -\tau_{is}) \right)}
  {\sum_{k=0}^J \exp\left(\sum_{s=1}^k \alpha_{is} (\theta_p - \tau_{is}) \right)}
\end{equation}
with
\begin{equation}
\sum_{s=1}^0 \alpha_{is} (\theta_p -\tau_{is}) := 0
\end{equation}

for notational convenience. The model results from the assumption that, given only two adjacent categories, the probability for the higher of the two is given by the dichotomous 2PL model:
\begin{equation}
P(Y = j; \theta_p, \boldsymbol{\tau_i}, \boldsymbol{\alpha_i} | Y_{pi} \in \{j, j-1\}) = \frac{\exp(\alpha_{ij} (\theta_p -\tau_{ij}))}{1 + \exp(\alpha_{ij}  (\theta_p -\tau_{ij}))}.
\end{equation}

The parameter $\theta_p$ denotes the ability of the $p$th person. The higher $\theta_p$, the higher the probability for reaching higher categories \cite{agresti2010}. The parameter vector $\boldsymbol{\tau_i} = (\tau_{i1}, ..., \tau_{iJ})$ denotes the so called \emph{thresholds} of the $i$th item. Thresholds in the GPCM can be interpreted as follows. If a person has an ability equal to $\tau_{ij}$, the probabilities for this person of achieving category $j$ and $j-1$ are equal, that is $\pi_{pij} = \pi_{pi(j-1)}$ (i.e. $\tau_{ij}$ are the intersection points of the respective category response curves). The higher the thresholds, the lower the probability for achieving higher categories. This does not imply, however, that thresholds must be ordered in the sense that $\tau_1 \leq \tau_2 \leq ... \leq \tau_J$ \cite{adams2012}. Rather, all thresholds may take on any real value. The vector $\boldsymbol{\alpha_i} = (\alpha_{i1}, ..., \alpha_{iJ})$, where $a_{ij} > 0$, denotes the positive \emph{discrimination} parameters of the $i$th item. The higher $\alpha_{ij}$, the steeper the curve favoring the $j$th category over the $(j-1)$th category with increasing $\theta_p$. Note that the GPCM is sometimes defined as having only a single discrimination parameter, which is assumed constant across categories, instead of a vector of discriminations. As the latter is a generalization of the former, all of our results may be applied to GPCMs with a single discrimination parameter, as well. 

Since the focus of the present paper is on optimal test designs, we consider the item parameters $\boldsymbol{\tau_i}$ and $\boldsymbol{\alpha_i}$ to be known (i.e. chosen by the design) and the person parameters $\theta_p$ as the quantities to be estimated. This situation arises in calibrated tests where item parameters were estimated in prior studies and the design question is to select the most appropriate items from the item pool. Assuming conditional independence of the items for given $\theta_p$, the joint density of all items is simply the product of the single-item densities (see \citeNP{vanderlinden2010b} for a detailed discussion about conditional independence in IRT models). For notational convenience, we will suppress the indices $p$ and $i$ in the following where appropriate. 

The original PCM was first derived by \citeA{rasch1961} and subsequently by \citeA{andersen1973}, \citeA{andrich1978a}, \citeA{masters1982}, and \citeA{Fischer1995} each with a different but equivalent formulation (considering the special case of discrimination parameters fixed to one; c.f. \citeNP{Fischer1995,adams2012}). \citeA{andersen1973} and \citeA{Fischer1995} derived the PCM in an effort to find a model that allows the independent estimation of person and item  parameters -- a highly desirable property -- for ordinal variables. \citeA{andrich1978a,andrich2005} provided another derivation: When two dichotomous processes are independent, four results can occur: $(0,0), (1,0), (0,1), (1,1)$. Using the Rasch model for each of the two processes, the probability of the combined outcomes is given by the Polytomous Rasch Model (\citeNP{andersen1973}; c.f. \citeNP{wilson1992,wilson1993}). When thinking of these processes as steps between ordered categories, $(0,0)$ corresponds to $Y = 0$, $(1,0)$ corresponds to $Y=1$, and $(1,1)$ corresponds to $Y = 2$. The event $(0,1)$, however, is assumed to be impossible because the second step cannot be successful when the first step was not. For an arbitrary number of ordered categories, \citeA{andrich1978a} proved that the Polytoumous Rasch Model becomes the PCM (with discrimination parameters fixed to one) when considering the set of possible events only. The GPCM, which generalizes the 2PL model to more than two ordinal response categories, was later proposed by \citeA{muraki1992}.


\section{3. Locally Optimal Designs}

In the context of optimal test designs for the GPCM, an experimental design is the set $\xi := \{ (\boldsymbol{\tau_1}, \boldsymbol{\alpha_1}), (\boldsymbol{\tau_2}, \boldsymbol{\alpha_2}), ..., (\boldsymbol{\tau_N}, \boldsymbol{\alpha_N}) \}$ of known parameters of the $N$ administrated items. It can be chosen by the experimenter in order to optimize the information of the experiment. Note that each $\boldsymbol{\tau_i}$ and each $\boldsymbol{\alpha_i}$ ($i \in \{1, ..., N\}$) is a vector of length $J$. The set of all designs is denoted as $\Xi$. When only considering a single item, we simply write $(\boldsymbol{\tau}, \boldsymbol{\alpha})$ instead of $\xi$ and drop the index $i$.

There are several optimal design criteria discussed in the literature (see \citeNP{atkinson2007} for an overview), but arguably the most common criterion is D-optimality. A design $\xi$ is called D-optimal if it minimizes the determinant of the estimator's covariance matrix. In the present context, the parameter $\theta$ to be estimated is the unidimensional person parameter so that the determinant of the covariance matrix reduces to the unidimensional variance. Thus, we can define a D-optimal design as
\begin{equation}
\xi^* := \argmin_{\xi \in \Xi} \Var (\hat{\theta}; \xi),
\end{equation}
where $\hat{\theta}$ is an estimator of the parameter $\theta$. According to the Cram{\'e}r-Rao bound \cite{cramer2016, rao1992}, the variance of an unbiased estimator cannot be smaller than the inverse of the Fisher-Information $M$: 
\begin{equation}
\Var (\hat{\theta}, \xi) \geq M(\theta, \xi)^{-1}
\end{equation}
Thus, for (asymptotically) efficient estimators, we can equivalently maximize the Fisher-Information instead of minimizing the estimator's variance. Under certain regularity conditions (cf. \citeNP{lehmann_casella2006}), $M$ can be written as
\begin{equation}
\label{Mdef}
M(\theta, \xi) := \Var \left[ \frac{\dd}{\dd \theta} \log f(Y; \theta, \xi)\right] = - \E \left[\frac{\dd^2}{\dd \theta^2} \log f(Y; \theta, \xi) \right],
\end{equation}
where $f(Y; \theta, \xi)$ is the product density of all items given the ability parameter $\theta$ and the design $\xi$ assuming conditional independence of the items. Due to additivity of the Fisher-Information in this case, we could equivalently write (\ref{Mdef}) as the sum of the single item Fisher-Informations. Most often in optimal design theory, the aim is to minimize the variance of the estimator. In case of asymptotic efficient estimators (such as maximum-likelihood estimators), which by definition meet the Cram{\'e}r-Rao bound, this is equivalent to maximizing the Fisher-Information and we will use this approach in the derivation of optimal designs in the present paper.

\begin{proposition}
\label{MPCM}
Defining $A_j := \sum_{s=1}^j \alpha_s$, the Fisher information $M(\theta, \boldsymbol{\tau}, \boldsymbol{\alpha})$ of a single item following the GPCM is given by
\begin{equation}
M(\theta, \boldsymbol{\tau}, \boldsymbol{\alpha}) = \sum_{j = 1}^{J} A_j^2 \pi_j - \left( \sum_{j = 1}^{J} A_j \pi_j \right)^2.
\end{equation}
\end{proposition}

If we treat $M$ as a function of the probability vector $\boldsymbol{\pi} = (\pi_0, ..., \pi_J)$, we can obtain the following theorem as an immediate result of Proposition \ref{MPCM}.

\begin{theorem}
\label{Mmax}
The Fisher information $M$ of a single item treated as a function of the probabilities $\boldsymbol{\pi}$ has a unique global maximum in $\{ \boldsymbol{\pi} \in \mathbb{R}_+^{J+1} \; | \; \sum_{j = 0}^J \pi_j = 1\}$, which is given by $\pi^*_0 = \pi^*_J = 1/2$ and $\pi^*_1 = \pi^*_2 = ... = \pi^*_{J-1} = 0$ for given $\boldsymbol{A} = (A_1, ..., A_J)$.
\end{theorem}

Theorem \ref{Mmax} implies that an item following the GPCM would be optimal for a person with ability $\theta_0$, when the vectors $\boldsymbol{\tau}$ and $\boldsymbol{\alpha}$ are chosen so that only the first and last category have non-zero probability and are equally likely. In other words, such an item is optimal if it is from the dichotomous 2PL model (with scalar discrimination $\alpha = A_J$). Due to additivity of the Fisher information, all items of a locally optimal design $\xi^*$ for the PCM have to satisfy the above condition.

We want to briefly illustrate the efficiency gain achieved by optimal items with a simple example. Suppose we administer a four category item with threshold vector $\boldsymbol{\tau} = (-1, 0, 1)$ and discrimination vector $\boldsymbol{\alpha} = (1, 1, 1)$. Then, participants with low ($\theta = -1$), medium ($\theta = 0$), and high ($\theta = 1$) ability achieve the categories 0 to 3 with probabilities $\pi_l \approx (0.41, 0.41, 0.15, 0.02)$, $\pi_m \approx (0.13, 0.37, 0.37, 0.13)$, and $\pi_h \approx (0.02, 0.15, 0.41, 0.41)$, respectively.  Such values are not too uncommon in practice. When compared to an item with optimal thresholds and $\boldsymbol{\alpha} = (1, 1, 1)$, the ratio of the Fisher-information is $M^* / M_m \approx 2.86$ for participants with medium ability, and $M^* / M_l = M^* / M_h \approx 3.75$ for participants with low or high ability. Hence, the variance of an efficient estimator is about three to four times higher when applying this item rather than optimal items tailored to the participants' abilities.

Unfortunately, there exists no item with a threshold vector $\boldsymbol{\tau^*}$ that leads to the optimal probabilities of Theorem \ref{Mmax}, but it may be approximated:

\begin{proposition}
\label{tau_star_properties}
Define $\boldsymbol{\tilde{\tau}_\alpha}(c) := (\alpha_J c, 0, 0, ..., -\alpha_1 c) \in \mathbb{R}^J$ with $c > 0$ and 
\begin{equation}
C^J_\alpha := \bigg\{x \in \mathbb{R}^J \; | \; \forall k \in \{1, ..., J\} : \alpha_{k} x_k = - \alpha_{J-k+1} x_{J-k+1} \bigg\}
\end{equation}
for some discrimination vector $\boldsymbol{\alpha} \in \mathbb{R}^J_+$.

\begin{enumerate}[label=(\alph*)]
\item It is $\boldsymbol{\tilde{\tau}_\alpha}(c) \in C^J_\alpha$.

\item For each $\theta \in \mathbb{R}$ there is some $\rho = \rho(\theta, \boldsymbol{\alpha})$ with $0 < \rho < 1$ such that 
$\pi_J(\theta, \boldsymbol{\tilde{\tau}_\alpha}(c), \boldsymbol{\alpha}) \rightarrow \rho$, $\pi_0(\theta, \boldsymbol{\tilde{\tau}_\alpha}(c), \boldsymbol{\alpha}) \rightarrow 1 - \rho$ and $\pi_j(\theta, \boldsymbol{\tilde{\tau}_\alpha}(c), \boldsymbol{\alpha}) \rightarrow 0$ for $0 < j < J$ as $c \rightarrow \infty$.

\item It is $\rho(0, \boldsymbol{\alpha}) = 1/2$.

\item For $\tau \in C^J_\alpha$ with symmetric $\boldsymbol{\alpha}$ satisfying $\alpha_k = \alpha_{J-k+1}$ it holds that $\pi_j(\theta, \boldsymbol{\tau}, \boldsymbol{\alpha}) = \pi_{J-j}(-\theta, \boldsymbol{\tau}, \boldsymbol{\alpha})$ for all $j \in \{0, ..., J \}$.
\end{enumerate}
\end{proposition}

According to Proposition \ref{tau_star_properties} (b) and (c), we can formally define $\boldsymbol{\tau^*}$ as $\boldsymbol{\tau^*} := \lim_{c \rightarrow \infty} \boldsymbol{\tilde{\tau}_\alpha}(c)$ when assuming $\theta_0 = 0$, which can be achieved by a location shift without loss of generality. Proposition \ref{tau_star_properties} may be generalized to arbitrary $\theta_0$, but because of the location shift argument, this is not required for the present paper. The sequence $\boldsymbol{\tilde{\tau}_\alpha}(c)$ is not the only one that satisfies the criteria (a) -- (d) of Proposition \ref{tau_star_properties}. However, for the purpose of the present paper, it is completely sufficient that we know at least one such sequence exists. Of course, it is impossible to create items with a threshold vector exactly equal to $\boldsymbol{\tau^*}$. Such an optimal item would have infinite and negative infinite first and last thresholds, respectively, while all other thresholds would not be estimable. In other words, for a fixed number of more than two response categories, a locally optimal design does not exist. However, one may choose reasonably large values of $c$ so that the GPCM closely approximates the dichotomous 2PL model with corresponding locally optimal design. The latter is well known as a one-point design, in which all items have difficulty equal to the persons ability. If additionally the discrimination is treated as part of the design, items should optimally have as high as possible discrimination. This may also be inferred directly from Proposition \ref{MPCM} when considering the special case of only two response categories.

\section{4. Bayesian Optimal Designs}

From a Bayesian perspective, a locally optimal design is the result of setting all prior mass to one point $\theta_0$. This seems quite infeasible and is only justified if the true ability is not too far away from $\theta_0$ or otherwise the locally optimal design in $\theta_0$ might perform poorly. One solution is to account for the a-priori uncertainty in the ability parameter by imposing a non-degenerate weight distribution $\Pi$ over the parameter space $\Theta$. An optimal design $\xi^* \in \Xi$ taking this uncertainty into account is called a Bayesian optimal design. It is obtained via some sort of averaging over the locally optimal designs' Fisher information, where the Fisher information of each locally optimal design is weighted according to the weight distribution \cite{atkinson2007}. In the present paper, we will consider two common Bayesian design criteria that each performs a different kind of averaging \cite{firth1997}:
\begin{align}
\psi_0(\xi) &:= \E\left[\log(M(\theta, \xi))\right], \\
\psi_{-1}(\xi) &:= - \E\left[M(\theta, \xi)^{-1}\right],
\end{align}
where 
\begin{equation}
\E\left[X\right] = \int X \, \Pi(\theta) \, \dd \theta
\end{equation}
is the expected value of a random variable $X$ under the distribution $\Pi$ of $\theta$. The criterion $\psi_0$ considers the average logarithm of the Fisher information, while the criterion $\psi_{-1}$ aims at minimizing the average asymptotic variance. Of those two, $\psi_0$ is mathematically more convenient and has a natural Bayesian interpretation \cite{chaloner1989}. Both criteria are concave with respect to the experimental design $\xi$ \cite{firth1997}. Note that, in general, these design criteria are not fully Bayesian in the sense that they do not require computing a posterior distribution. Instead, they borrow the concept of prior distributions to obtain robust designs for frequentist analysis via maximum likelihood methods. Hence, the here discussed designs are sometimes referred to as 'pseudo-Bayesian' designs \cite{firth1997}. For a review of fully Bayesian design criteria see \citeA{chaloner1995}.

In theory, the weight distribution $\Pi$ of $\theta$ may have any form. In practice, however, one typically uses symmetric weights, both for mathematical convenience and because there is rarely any prior information available that favors an asymmetric weight distribution for an unbounded parameter such as the ability $\theta$. Accordingly, we will focus on symmetric weight distributions in the following. Within this class, there are -- among others -- the continuous uniform, normal, and logistic distribution, each forming a so called location scale family of distributions. 

For the Rasch model, Bayesian optimal designs have already been investigated in \citeA{grasshoff2012} and the present section aims at generalizing these results to the GPCM. The following Lemma underlines the relevance of locally optimal designs in the context of Bayesian optimal design criteria.

\begin{lemma}
\label{BOD_one_point}
If the weight distribution $\Pi$ is symmetric around some ability $\theta_0$, $\boldsymbol{\alpha}$ is fixed to any vector in $\mathbb{R}_+^J$, and only $\boldsymbol{\tau}$ is considered variable, the Bayes optimal one-point design with respect to $\psi_0$ and $\psi_{-1}$ is the locally optimal design for $\theta_0$ that is $\boldsymbol{\tau} = \boldsymbol{\tau^*}$.
\end{lemma}

According to Lemma~\ref{BOD_one_point}, the same problem that occurs for the locally optimal design occurs for the optimal one-point design under the Bayesian criteria as well. That is, it does not exist for a fixed number of more than two response categories and when considering the number of categories variable, the GPCM reduces to the dichotomous 2PL model. Accordingly, we need to investigate the optimal discrimination parameter $\alpha$ in case of only two response categories:

\begin{lemma}
\label{BOD_one_point_alpha}
If the weight distribution $\Pi$ is symmetric around some ability $\theta_0$, $\boldsymbol{\tau}$ is fixed to $\boldsymbol{\tau^*}$, and $\alpha \in \mathbb{R}_+$ is considered variable, it holds that:
\begin{enumerate}[label=(\alph*)]
\item $\psi_0$ is maximized by the unique solution $\alpha^*_0$ of
\begin{equation}
\alpha \int \theta \ \pi_1(\theta, \boldsymbol{\tau^*}, \alpha) \, \Pi(\theta) \, \dd \theta = 1.
\end{equation}
\item $\psi_{-1}$ is maximized by the unique solution $\alpha^*_{-1}$ of
\begin{equation}
\alpha \int \frac{\theta}{\pi_0(\theta, \boldsymbol{\tau^*}, \alpha)} \, \Pi(\theta) \, \dd \theta = \int \frac{1}{\pi_1(\theta, \boldsymbol{\tau^*}, \alpha) \ \pi_0(\theta, \boldsymbol{\tau^*}, \alpha)} \, \Pi(\theta) \, \dd \theta.
\end{equation}
\end{enumerate}
\end{lemma}

Lemma \ref{BOD_one_point_alpha} is in line with our intuition that the optimal $\alpha$ decreases with increasing scale of $\Pi$.  Lemmas \ref{BOD_one_point} and \ref{BOD_one_point_alpha} do not state under which conditions the optimal one-point design is Bayes optimal. To answer this question, we have to make use of approximate design theory (cf. \citeNP{kiefer1974}) and introduce so called sensitivity functions measuring the quality of a design (cf. \citeNP{grasshoff2012}). If the sensitivity function is uniformly bounded at a certain value for a given design $\xi^*$, no improvement is possible and $\xi^*$ will be optimal. For non-linear models, sensitivity functions of D-optimal designs are given in \citeA{firth1997}. In case of the GPCM, these sensitivity functions can be written as
\begin{equation}
\phi_v(\boldsymbol{\tau}, \boldsymbol{\alpha}, \xi) = \frac{E\left[M(\xi, \theta)^{v} M(\xi, \theta)^{-1} M(\boldsymbol{\tau}, \boldsymbol{\alpha}, \theta) \right]}{E\left[M(\xi, \theta)^{v}\right]} - 1,
\end{equation}
where $v \in \{-1, 0\}$ corresponds to the two Bayesian design criteria discussed above. Since, for the GPCM, the design criteria are concave, Bayes optimality of a design $\xi^*$ is equivalent to 
\begin{equation}
\label{CET}
\sup_{(\boldsymbol{\tau}, \boldsymbol{\alpha}) \in \mathcal{X}}\phi_v(\boldsymbol{\tau}, \boldsymbol{\alpha}, \xi^*) = 0,
\end{equation}
where $\mathcal{X}$ denotes the design region of a single item. Moreover, when $\xi^*$ is optimal, $\phi_v(\boldsymbol{\tau}, \boldsymbol{\alpha}, \xi^*) = 0$ holds if and only if $(\boldsymbol{\tau}, \boldsymbol{\alpha})$ is a design point of $\xi^*$. It may now be asked, under which conditions the optimal one-point design is also Bayes optimal. 

\begin{theorem}
\label{BOD_PCM}
Necessary conditions for the threshold $\boldsymbol{\tau^*}$ to be Bayes optimal under a symmetric weight distribution $\Pi$ and fixed discrimination $\alpha$ are 
\begin{equation}
\label{ncond0}
\psi_0: \quad \int \pi_J(\theta, \boldsymbol{\tau^*}, \boldsymbol{\alpha}) \, \pi_0(\theta, \boldsymbol{\tau^*}, \boldsymbol{\alpha}) \, \Pi(\theta) \, \dd \theta \geq \frac{1}{6}
\end{equation}
\begin{equation}
\label{ncondM1}
\psi_{-1}: \quad \int (\pi_J(\theta, \boldsymbol{\tau^*}, \boldsymbol{\alpha}) \, \pi_0(\theta, \boldsymbol{\tau^*}, \boldsymbol{\alpha}))^{-1} \, \Pi(\theta) \, \dd \theta \leq 6
\end{equation}
\end{theorem}

For a fixed discrimination of $\boldsymbol{\alpha} = 1$, Theorem \ref{BOD_PCM} provides the same conditions for the GPCM as can be found for the Rasch model (see \citeNP{grasshoff2012}). By Jensen's inequality, the necessary condition for $\psi_{-1}$ implies the condition for $\psi_0$, that is the former condition is more restrictive. For scale families of symmetric weight distributions, Theorem \ref{BOD_PCM} may be further extended. 

\begin{theorem}
\label{BOD_PCM_symmetric}
Let $(\Pi_s)_{s > 0}$ be a scale family of symmetric distributions around some ability $\theta_0$.
\begin{enumerate}[label=(\alph*)]
\item There exist unique values $s_{-1}$ and $s_0$ with $0 < s_{-1} \leq s_0 \leq \infty$, such that the necessary conditions of Theorem \ref{BOD_PCM} on $\psi_{-1}$ and $\psi_0$ are satisfied for all $s \leq s_{-1}$ and $s \leq s_0$, respectively, and violated otherwise.
\item There exists a unique threshold $\tilde{s}_0$, $0 \leq \tilde{s}_0 < s_0$, such that the optimal one-point design $\tau^*$ is Bayes optimal with respect to $\psi_0$ for all weight distributions $\Pi$ with $s \leq \tilde{s}_0$ and fails to be optimal for $s > \tilde{s}_0$.
\end{enumerate}
\end{theorem}

According to Theorem \ref{BOD_PCM_symmetric}, the optimal one-point design for a symmetric scale distribution remains Bayes optimal as long as the scale parameter does not exceed a certain value that depends on the design criterion, the (now scalar) discrimination $\alpha$, and on the location scale family. As $\int \pi_J \pi_0 \, \Pi(\theta) \, \dd \theta$ decreases monotonically with increasing $\alpha$, the values $\tilde{s}_0$, $s_0$, and $s_{-1}$ decrease monotonically with increasing $\alpha$ as well. For some common distribution families, examples are discussed in more detail in the upcoming section. 

Lastly, we want to turn our attention to Bayes optimal designs for the dichotomous 2PL model, when both $\tau$ and $\alpha$ are allowed to vary. It is of particular interest, in which cases the optimal one-point designs derived in Lemma \ref{BOD_one_point_alpha} are also Bayes optimal. In next section, we show numerically that for certain important weight distributions under $\phi_0$ or $\phi_{-1}$, the optimal one-point design is never Bayes optimal. This implies that the Bayes optimal designs under $\phi_0$ or $\phi_{-1}$ (if existent) have to have at least two distinct design points. Whether this holds in general for all symmetric weight distributions requires further investigation.


\section{5. Examples}

Among others, \citeA{grasshoff2012} discuss Bayes optimal designs for the Rasch model in case of the continuous uniform, normal, and logistic families. They also derive the family specific values of $s_{-1}$ and $s_0$ (named $\tau_1$ and $\tau_2$, respectively, in their paper) for the Rasch model (i.e. for discrimination fixed to $1$). We have already established that, when using the optimal one-point design, the GPCM is equivalent to the 2PL model. We can obtain $s_{0}$ and $s_{-1}$ for the GPCM by dividing the respective values for Rasch model by the sum $A_J$ of all discrimination parameters. Theoretically, this allows to numerically investigate, whether the necessary conditions on the optimality of the one-point design are also sufficient, by computing the sensitivity function across the design space. Practically, this becomes increasingly complicated for items with more response categories, since the design space is of dimension $J$ (or $2J$ when also varying $\boldsymbol{\alpha}$). Accordingly, we will focus on items with up to three response categories, which may also be visualized conveniently. For each example, we will consider two cases: (a) three response categories with varying $\boldsymbol{\tau}$ and discrimination parameters fixed to one (hence $A_j = j$ for all $j$) as well as (b) two response categories with varying $\tau$ and $\alpha$. \\

\begin{figure}[htbp]
	\centering
  \includegraphics[width=0.99\textwidth]{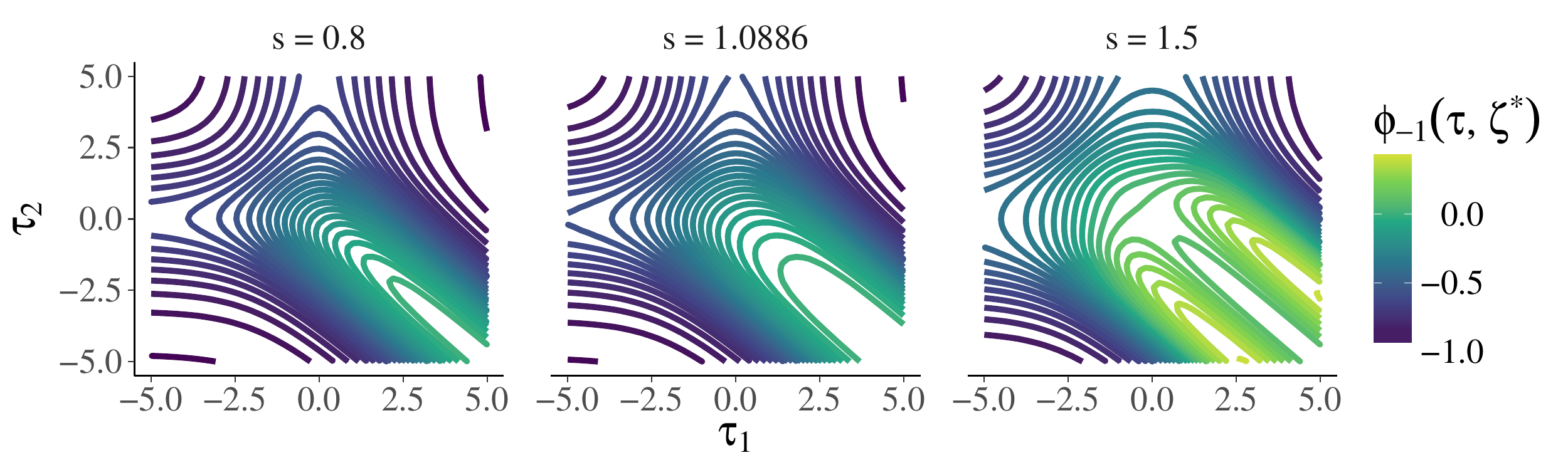}
	\caption{Contour of the sensitivity function for $\psi_{-1}$ in case of three response categories for a uniform weight distribution with varying values of the scale parameter $s$.}
	\label{unif_m1}
\end{figure}

\begin{figure}[htbp]
	\centering
  \includegraphics[width=0.99\textwidth]{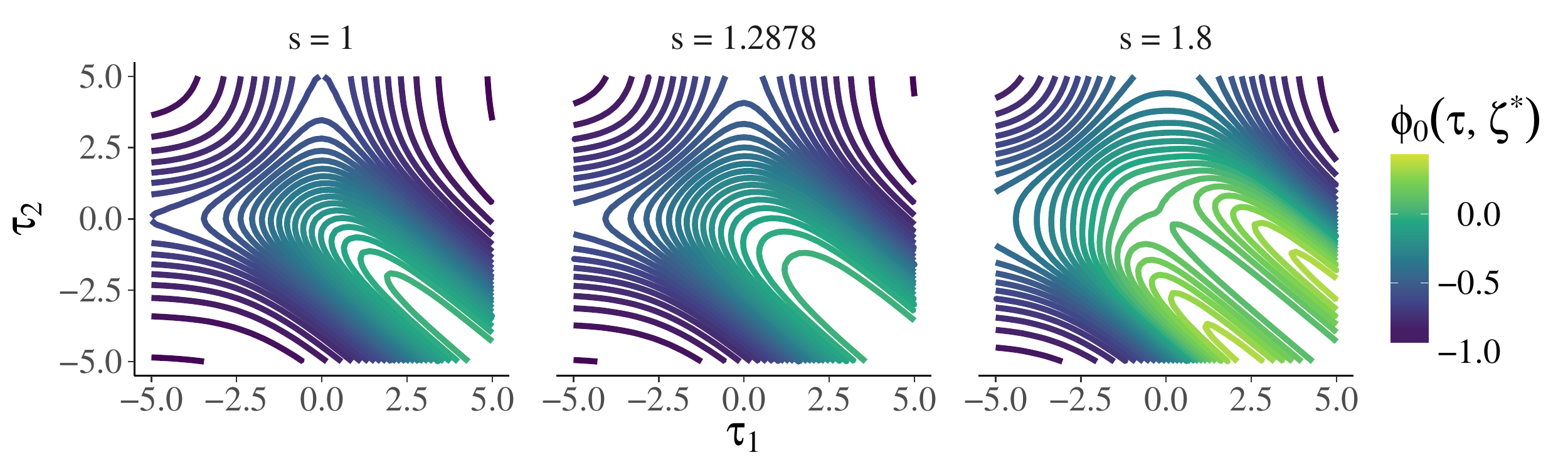}
	\caption{Contour of the sensitivity function for $\psi_{0}$  in case of three response categories for a uniform weight distribution with varying values of the scale parameter $s$.}
	\label{unif_0}
\end{figure}

\emph{Example 1} (Continuous uniform distribution): For the continuous uniform distribution with support in $[\theta_1, \theta_2] = [\theta_0 - s, \theta_0 + s]$ and constant density $f(\theta) = \frac{1}{2s}$, we have $s_{-1} \approx 2.1773 / J$ and $s_{0} \approx 2.5757 / J$. The scalars $2.1773$ and $2.5757$  are numeric solutions of the equations  $\exp(s) − \exp(-s) = 4s$ and $(3 - s) \exp(s) = 3 + s$, respectively \cite{grasshoff2012}. For $J = 2$, the sensitivity function is displayed in Figure \ref{unif_m1} and \ref{unif_0}, with $s_{-1}$ and $s_{0}$ being displayed in the center. In the Figures, yellowish regions are values of the sensitivity function that are greater than zero hence indicating non-optimality of the one-point design. At $s_{-1}$ and $s_{0}$, or smaller values, the sensitivity functions do not exceed zero, and are unimodal with maximal value zero at $\boldsymbol{\tau^*}$. This demonstrates, at least numerically and for $J = 2$, that the necessary conditions of Theorem \ref{BOD_PCM} are also sufficient for the continuous uniform weight distribution. When $s$ is further increased, the one-point design becomes visibly non-optimal and a two-point design appears to be favored. 

\begin{figure}[htbp]
	\centering
  \includegraphics[width=0.99\textwidth]{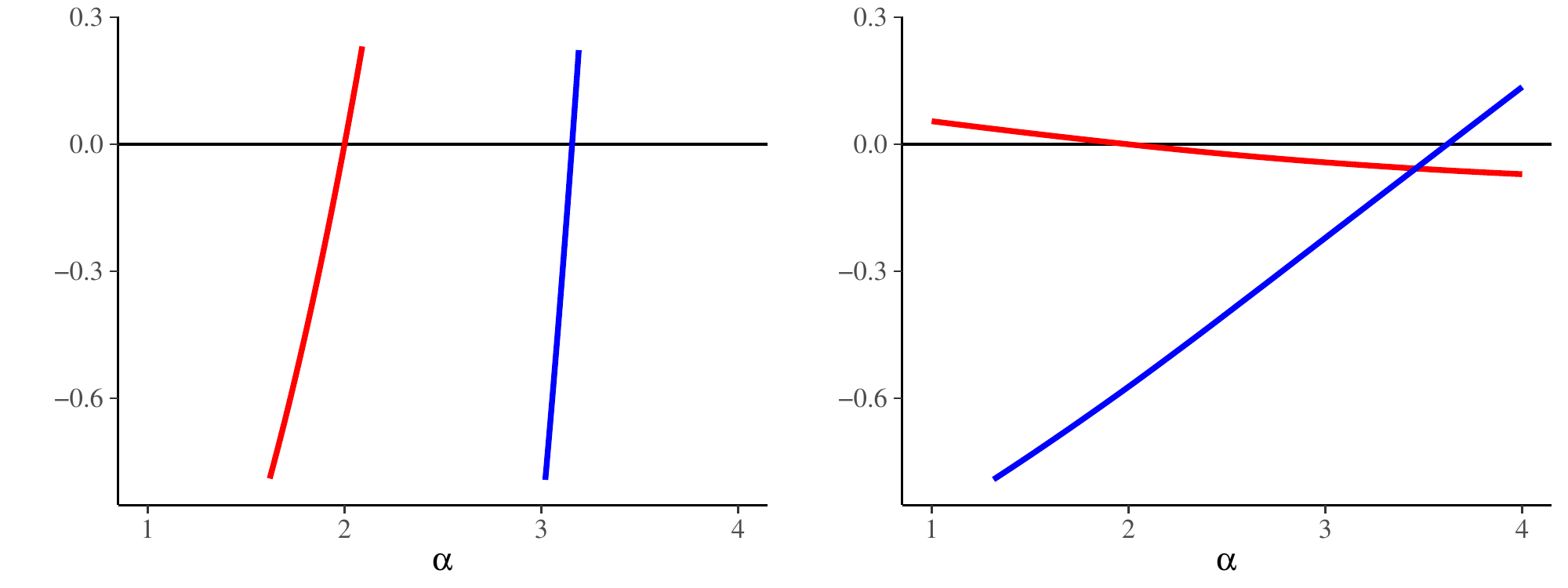}
	\caption{Illustration of the optimal one-point discrimination parameter of the 2PL model in relation to the necessary condition for Bayes optimality. Left: criterion $\phi_{-1}$ with a uniform weight distribution with scale $s = 2.1773$. Right: criterion $\phi_{0}$ with a uniform weight distribution with scale $s = 2.5757$. The intersection of the blue and black line indicate $\alpha^*_{-1}$ and $\alpha^*_{0}$, respectively, while the intersection of the red and black line indicate $\alpha^+_{-1}$ and $\alpha^+_{0}$, respectively. Note that $\alpha^*_{-1} > \alpha^+_{-1}$ and $\alpha^*_{0} > \alpha^+_{0}$.}. 
	\label{optimal_alpha_unif}
\end{figure}

For the dichotomous 2PL model with varying $\tau$ and $\alpha$,  we have $\alpha^*_{-1} \approx 3.1560$ for $s = s_{-1} \approx 2.1773$. However, the necessary condition (\ref{ncondM1}) for $\phi_{-1}$ is only satisfied for $\alpha \leq \alpha^+_{-1} \approx 2$, implying that the optimal one-point design $(\tau^*, \alpha^*_{-1})$ is not Bayes optimal (see also Figure \ref{optimal_alpha_unif} left-hand side). The same result may be obtained for other values of the scale $s$ in the sense that $\alpha^*_{-1}$ exceeds $\alpha^+_{-1}$ independently of $s$. Similarly, we have $\alpha^*_{0} \approx 3.6186$ for $s = s_0 \approx 2.5757$, but the necessary condition (\ref{ncond0}) for $\phi_0$ is only satisfied for $\alpha \leq \alpha^+_0 \approx 2$, (see also Figure \ref{optimal_alpha_unif} right-hand side). Again, $\alpha^*_{0} > \alpha^+_{0}$ can be shown to hold for other values of $s$ as well.\\

\begin{figure}[htbp]
	\centering
  \includegraphics[width=0.99\textwidth]{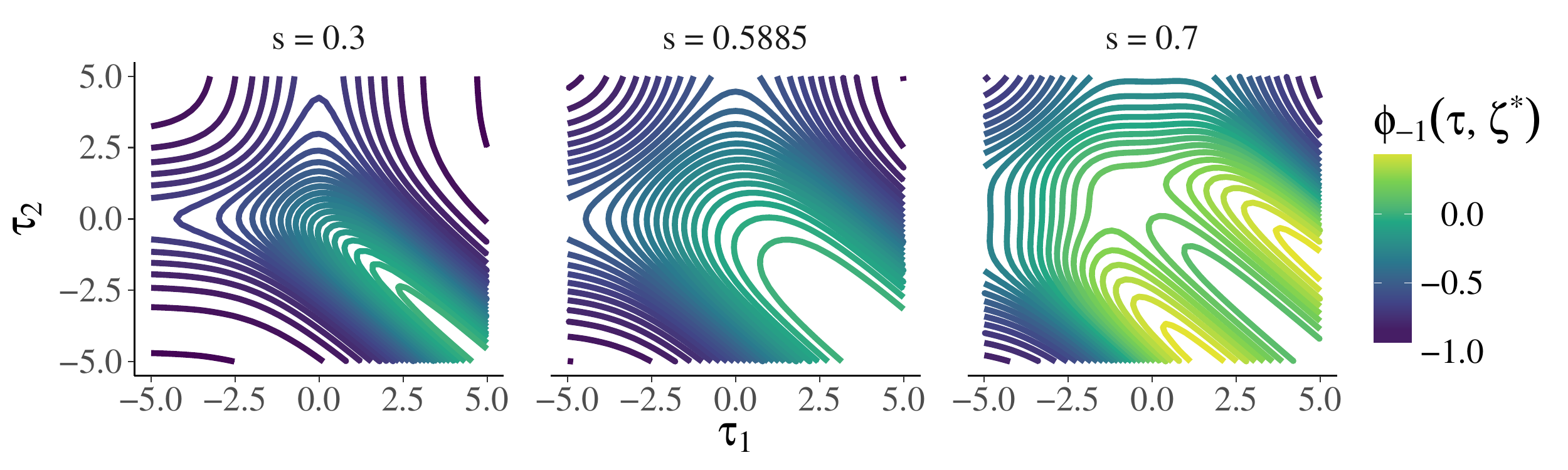}
	\caption{Contour of the sensitivity function for $\psi_{-1}$ in case of three response categories for a normal weight distribution with varying values of the scale parameter $s$.}
	\label{norm_m1}
\end{figure}

\begin{figure}[htbp]
	\centering
  \includegraphics[width=0.99\textwidth]{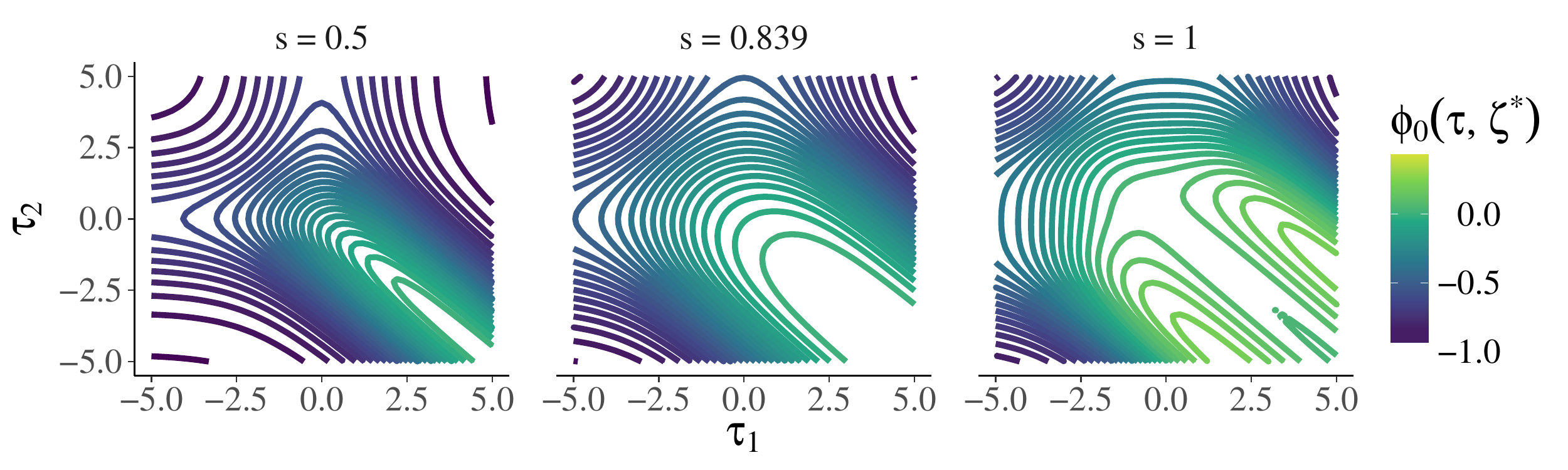}
	\caption{Contour of the sensitivity function for $\psi_{0}$ in case of three response categories for a normal weight distribution with varying values of the scale parameter $s$.}
	\label{norm_0}
\end{figure}

\emph{Example 2} (Normal distribution): Similar results are obtained for the normal weight distribution.  Here, the scale $s$ is simply the standard deviation parameter of the normal distribution. The critical values are approximately $s_{-1} \approx 1.177 / J$ and $s_{0} \approx 1.683 / J$. Again, numerical computation shows that the sensitivity functions do not exceed zero for $s \leq s_{-1}$ or $s \leq s_{0}$ (see Figure \ref{norm_0} and \ref{norm_m1}). Thus, the necessary conditions appear to be sufficient for the normal weight distribution when $J = 2$. Accordingly, the optimal one-point design remains Bayes optimal for items with three response categories as long as the standard deviation does not exceed the critical values. 

For the dichotomous 2PL model with varying $\tau$ and $\alpha$,  we have $\alpha^*_{-1} \approx 1.3586$ for $s = s_{-1} \approx 1.177$, but the necessary condition (\ref{ncondM1}) for $\phi_{-1}$ is only satisfied for $\alpha \leq \alpha^+_{-1} \approx 1$. Similarly, we have $\alpha^*_{0} \approx 1.7350$ for $s = s_0 \approx 1.683$, but the necessary condition (\ref{ncond0}) for $\phi_0$ is only satisfied for $\alpha \leq \alpha^+_0 \approx 1.002$. This may be shown to hold for other values of $s$ as well. \\

\begin{figure}[htbp]
	\centering
  \includegraphics[width=0.99\textwidth]{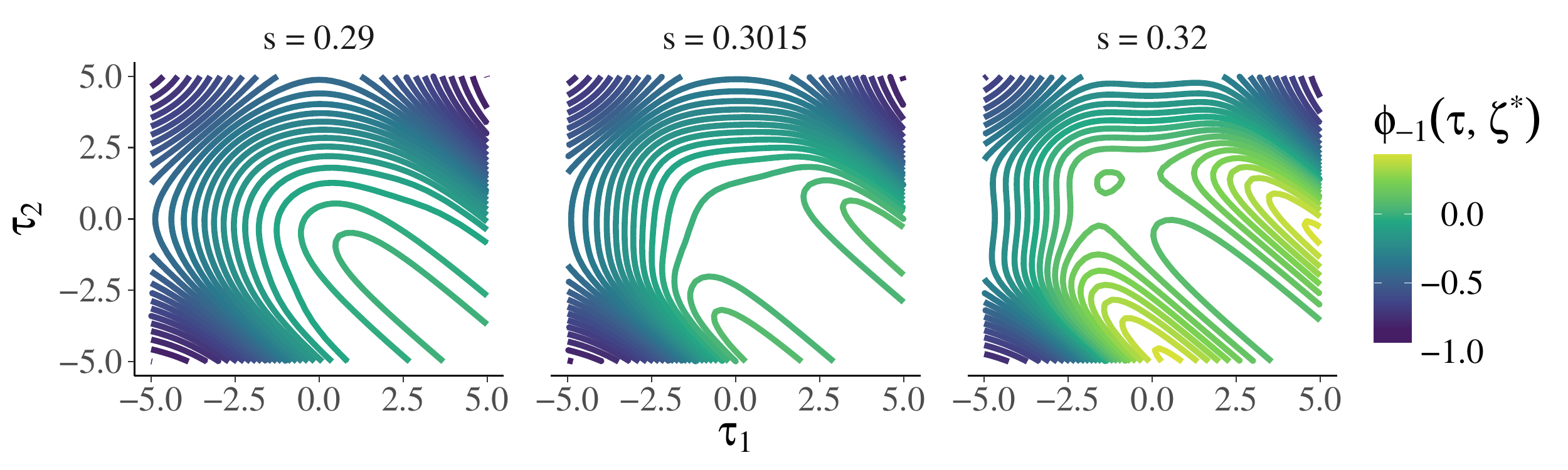}
	\caption{Contour of the sensitivity function for $\psi_{-1}$ in case of three response categories for a logistic weight distribution with varying values of the scale parameter $s$.}
	\label{logis_m1}
\end{figure}

\begin{figure}[htbp]
	\centering
  \includegraphics[width=0.99\textwidth]{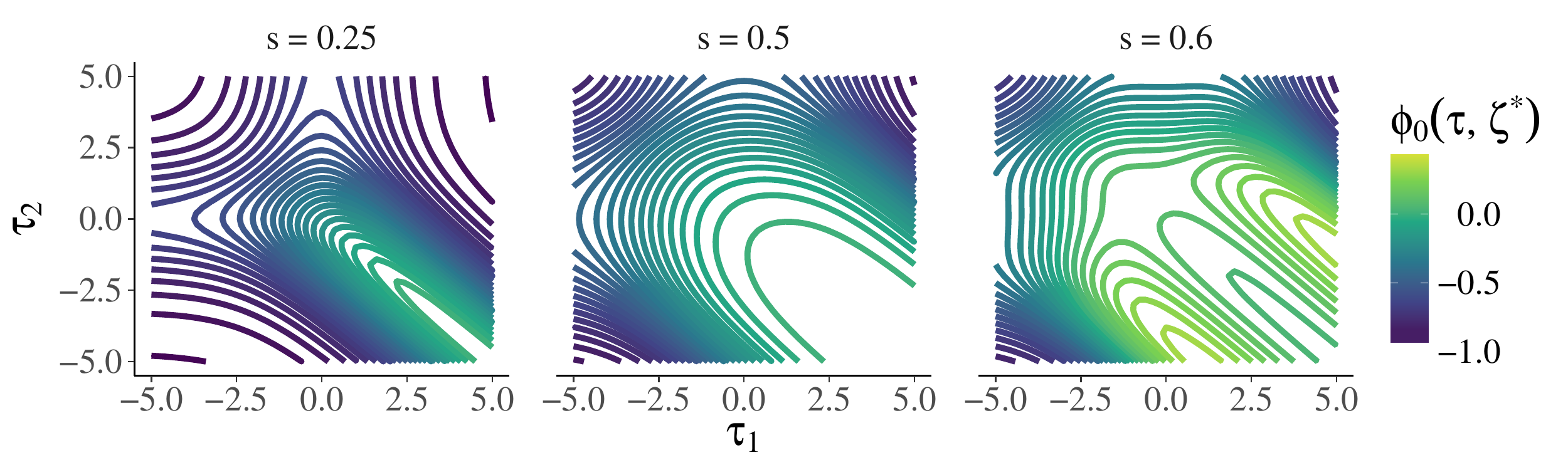}
	\caption{Contour of the sensitivity function for $\psi_{0}$ in case of three response categories for a logistic weight distribution with varying values of the scale parameter $s$.}
	\label{logis_0}
\end{figure}

\emph{Example 3} (Logistic distribution): For the logistic distribution with scale parameter $s$, the critical values are given by $s_{-1} \approx 0.603 / J$ and $s_0 = 1 / J$. According to Figure \ref{logis_0}, the necessary condition for $\psi_0$ is also sufficient if $J = 2$. However, the sensitivity function of $\psi_{-1}$ behaves somewhat differently for the logistic weight distribution in the sense that the necessary condition is apparently not sufficient (see Figure \ref{logis_m1} center). However, for scale values only slightly below $s_{-1}$, the optimal one-point design appears to be Bayes optimal (see Figure \ref{logis_m1} left hand side), indicating the the best possible sufficient condition, it very close to the necessary condition derived in the present paper. 

For the dichotomous 2PL model with varying $\tau$ and $\alpha$,  we have $\alpha^*_{-1} \approx 2.6518$ for $s = s_{-1} \approx 0.603$, but the necessary condition (\ref{ncondM1}) for $\phi_{-1}$ is only satisfied for $\alpha \leq \alpha^+_{-1} \approx 1.953$. Similarly, we have $\alpha^*_{0} \approx 2.9217$ for $s = s_0 \approx 1$, but the necessary condition (\ref{ncond0}) for $\phi_0$ is only satisfied for $\alpha \leq \alpha^+_0 \approx 1.6868$. This may be shown to hold for other values of $s$ as well. \\

\section{6. Discussion}

In the present paper, we investigated locally optimal designs as well as Bayes optimal designs for the GPCM. We demonstrated that for a fixed number of more than two response categories, the locally optimal design does not exist. However, when allowing the number of categories to vary, the locally optimal design exists and assigns non-zero probability only to the first and last response category independently of a person's ability. That is, when using such a design, the GPCM reduces to the ordinary 2PL model. This perhaps paradoxical result will be discussed in more detail later on.


When fixing the discrimination parameters, the locally optimal one-point design proved to be relatively robust. In particular, when allowing the number of response categories to vary, it turned out to be Bayes optimal under the two most common Bayesian design criteria and symmetric weight distribution with not too large variation. It has to be emphasized that the maximal variation such that it remained Bayes optimal for both criteria can be considered as sufficiently large for application. For instance, when using a normal weight distribution with standard deviation of $s \approx 0.59$ as well as items with three response categories and a discrimination of one, subjects in the lowest $2.5\%$ ability-region (according to the weight distribution) completely solve the task only with probability $9\%$ or less, while subjects in the highest $2.5\%$ ability-region completely solve the task with probability $91\%$ or more. Thus, even for substantial variation in the subjects' abilities, one-point designs remain compatitively viable.

When focusing on dichotomous responses and thus on the 2PL model, we provided equations whose unique solutions are the optimal discrimination parameters in the set of one-point designs when applying symmetric weight distributions. However, for some important examples, we established numerically that these discrimination parameters do not lead to Bayes optimal designs. That is, any Bayes optimal design for the 2PL model with varying threshold and discrimination (if existent) has to have at minimum two distinct design points at least for the examples discussed in the present paper.

Despite being mathematically reasonable, our results reveal practical problems when trying to design items with optimal properties. An item with more than two response categories, on which subjects score on the first or last categories (nearly) always -- independent of their ability -- may be practically impossible to construct. One could only aim at reducing the probability of intermediate responses to approximate the optimal design. However, this would contradict the very idea of using ordinal items instead of dichotomous ones. The obvious solution to completely remove partial credit from the items and score only as true or false, might be tempting, but will definitely lead to information loss -- and thus less efficient designs -- for complicated items consisting of multiple tasks. This does not invalidate the GPCM in general as the lack of a reasonable locally optimal design does not affect other mathematical or psychometric properties of the model. Further, the Fisher information or other criteria may still be used to compare items or to choose items with the highest information among those being available. 

Still, it is of great interest to study optimal designs of alternative ordinal models such as the graded response model \cite{walker1967, vanderark2001} or the sequential model \cite{tutz1990, tutz2000}, which potentially favor items with properties much closer to real items. If so, one  should also evaluate their potential use in large scale educational assessment studies such as PISA, which had previously applied the GPCM. 

For future research, some open questions remain. It is still unclear how the optimal designs derived here perform for asymmetric weight distributions due to the complexity of the involved derivatives and integrals. Furthermore, it has to be noted that optimal test designs and optimal sampling designs are not equivalent as there are multiple threshold as well as discrimination parameters per item, but only one ability parameter per person. Accordingly, derivation of optimal sampling designs for the GPCM still remains an open topic. Finally, as discussed above, there are other reasonable ordinal models, whose optimal designs are potentially of great relevance in practice. While locally optimal designs for the graded response model have been derived in \citeA{schmidt2015}, Bayes optimal designs for this model as well as optimal designs in general for the sequential model remain to be investigated.

\vspace{\fill}\pagebreak

\section{7. Appendix}

The indexes $p$ and $i$ for persons and items are suppressed for notational convenience.

\begin{proof} \emph{(Proposition \ref{MPCM})}
We abbreviate the normalizing function in the denominator of (\ref{PCM}) as $z \, (= 1 / \pi_0)$ for convenience. For a single item, the response $Y \in \{0, ..., J\}$ is categorically distributed so that the log-likelihood $l$ equals:
\begin{align*}
l(Y; \theta, \boldsymbol{\tau}, \boldsymbol{\alpha}) := \log \left( f(Y; \theta, \tau, \alpha) \right) &= \sum_{j = 0}^{J} 1_{j}(Y) \log(\pi_j) + q \\
&= \sum_{j = 1}^{J} 1_{j}(Y) \log(\pi_j) + (1 - \sum_{j = 1}^{J} 1_{j}(Y)) \log(\pi_0) + q \\
&=  \sum_{j = 1}^{J} \left(1_{j}(Y) \log\left(\frac{\pi_j}{\pi_0}\right) \right) + \log(\pi_0) + q  \\
&= \sum_{j = 1}^{J} \left(1_{j}(Y) \sum_{s = 1}^j \alpha_s (\theta - \tau_s) \right) - \log(z) + q, \numberthis
\end{align*}
where $f$ is the density of $Y$ defined in (\ref{PCM}), $1_j$ is the indicator function for $Y = j$ and $q$ is a constant independent of $\theta$.  The first derivative of $l$ with respect to $\theta$ equals
\begin{align*}
\label{dldtheta}
\frac{\dd l(\theta, \boldsymbol{\tau}, \boldsymbol{\alpha})}{\dd \theta} &= \sum_{j = 1}^{J} A_j 1_{j}(Y) - \frac{\sum_{j = 1}^{J} A_j \exp \left(\sum_{s=1}^j \alpha_s (\theta - \tau_{s}) \right)}{z} \\
&= \sum_{j = 1}^{J} A_j (1_{j}(Y) - \pi_j)  \numberthis
\end{align*}
Thus, $M$ is equal to the variance of the random variable 
\begin{equation}
\label{Xdef}
X := \sum_{j = 1}^{J} A_J 1_{j}(Y).
\end{equation}
Since $(1_j(Y))^2 = 1_j(Y)$ and $1_j(Y) 1_k(Y) = 0$ for $j \neq k$, as well as $\E(1_j(Y)) = \pi_j$, we have
\begin{align*}
\label{VarX}
\Var(X) = \E(X^2) - \E(X)^2 = \sum_{j=1}^{J} A_j^2 \pi_j- \left(\sum_{j=1}^{J}A_j\pi_j \right)^2,
\end{align*}


which completes the proof. \hfill $\Box$
\end{proof}

\begin{proof} \emph{(Theorem \ref{Mmax})}
The random variable $X$ defined in (\ref{Xdef}), is bounded between $0$ and $A_J$. It is minimal if $Y = 0$ and maximal if $Y = J$. According to Popoviciu's inequality for variances \cite{popoviciu1965, peajcariaac1992}, $\Var(X)$ (and hence $M$) is maximal if and only if $\pi_0 = \pi_J = 1/2$ and $\pi_j = 0$ for $0 < j < J$.\hfill $\Box$
\end{proof}

\begin{proof} \emph{(Proposition \ref{tau_star_properties})} 

(a) This can easily be checked.

(b) We have
\begin{equation}
\sum_{k=0}^J \exp\left(\sum_{s=1}^k \alpha_s (\theta- \tilde{\tau}_{\alpha s}(c)) \right) =
1 + \left(\sum_{k=1}^{J-1} \exp\left(A_k \theta - \alpha_1 \alpha_J c) \right) \right) + \exp(A_J \theta) 
\end{equation}
Hence, for $c \rightarrow \infty$ it holds that $\pi_0 = 1 / (1 + \exp(A_J \theta))$, $\pi_J = \exp(A_J \theta) / (1 + \exp(A_J \theta))$ and $\pi_j = 0$  for all $1 \leq j < J$.

(c) This follows directly from the proof of (b).

(d) For $\tau \in C^J_\alpha$ with $\boldsymbol{\alpha}$ satisfying $\alpha_k = \alpha_{J-k+1}$ ($k \in \{1, ..., J\}$) we have
\begin{align*}
\exp \left(\sum_{s=1}^{J-j} \alpha_s (- \theta - \tau_{s}) \right) &= 
\exp \left(- A_{J - j} \theta - \sum_{s=1}^{J-j} \alpha_s \tau_{s} \right) \\
&= \exp \left(- (A_{J} - A_j) \theta - \sum_{s=1}^{j} \alpha_s \tau_{s} \right) \\
&= \exp(-A_J \theta) \exp \left(\sum_{s=1}^j \alpha_s ( \theta - \tau_{s}) \right) \numberthis
\end{align*}
since $A_J = A_j + A_{J-j}$ and $\sum_{s=1}^{J-j} \alpha_s \tau_{s} = \sum_{s=1}^{j} \alpha_s \tau_{s}$ for all $j \in \{0, ..., J\}$. Again for $\tau \in C^J_\alpha$ the denominator of $\pi_{J-j}(-\theta, \boldsymbol{\tau}, \boldsymbol{\alpha})$ can be written as
\begin{align*}
\sum_{k=0}^J \exp \left(\sum_{s=1}^{k} \alpha_s (- \theta - \tau_{s}) \right) 
&= \exp(-A_J \theta) \sum_{k=0}^J \exp \left( (A_{J} - A_k) \theta - \sum_{s=1}^{k} \alpha_s \tau_{s} \right) \\
&= \exp(-A_J \theta) \sum_{k=0}^J \exp \left( A_{J-k} \theta - \sum_{s=1}^{J-k} \alpha_s \tau_{s} \right) \\
&= \exp(- A_J \theta) \sum_{k=0}^J \exp \left( \sum_{s=1}^{k} \alpha_s (\theta -  \tau_{s}) \right). \numberthis
\end{align*}
Hence, we find
\begin{equation}
\pi_{J-j}(-\theta, \boldsymbol{\tau}, \boldsymbol{\alpha}) = \frac{\exp(-A_J \theta)}{\exp(-A_J \theta)} \, \pi_{j}(\theta, \boldsymbol{\tau}, \boldsymbol{\alpha}) = \pi_{j}(\theta, \boldsymbol{\tau}, \boldsymbol{\alpha})
\end{equation}
\hfill $\Box$
\end{proof}

In the following proofs, we will suppress the dependency of $M$, $\pi_j$, and $\Pi$ on $\theta$ for notational convenience. Whenever we apply the degenerate locally optimal design $\tau^*$, we implicitly reduce the number of categories to two, as in this case, only two response categories have positive probability. This way, the derived results remain valid without writing down any limits of the form $\lim_{\tau \rightarrow \tau^*}$ thus simplifying the notation.

\begin{proof} \emph{(Lemma \ref{BOD_one_point})}
For $j \geq i $ the derivatives of $\pi_j$ with respect to $\tau_i$ and $\alpha_i$ equal
\begin{align*}
\frac{\dd \pi_j}{\dd \tau_i} &= \alpha_i \left(\exp\left(\sum_{s=1}^j \alpha_s (\theta - \tau_s) \right) \left( \sum_{k = i}^J \exp\left(\sum_{s=1}^k \alpha_s (\theta - \tau_s) \right) - z \right) \right) / z^2 \\
&= \alpha_i \pi_j \left(\sum_{k = i}^J \pi_k - 1 \right) = 
\alpha_i \left( \pi_j \left(\sum_{k = i}^J \pi_k \right) - \pi_j \right). \numberthis
\end{align*}
For $j < i$ the derivatives equal
\begin{align*}
\frac{\dd \pi_j}{\dd \tau_i} &= \alpha_i \left(\exp\left(\sum_{s=1}^j \alpha_s (\theta - \tau_s) \right) \left( \sum_{k = i}^J \exp\left(\sum_{s=1}^k \alpha_s (\theta - \tau_s) \right) \right) \right) / z^2 \\
&= \alpha_i \pi_j \left(\sum_{k = i}^J \pi_k \right). \numberthis
\end{align*}

Thus, we find
\begin{equation}
\sum_{j = 1}^{J} \frac{\dd \pi_j}{\dd \tau_i} = \alpha_i \left( \left(\sum_{j=1}^{J} \pi_j \right) \left(\sum_{k = i}^J \pi_k \right) - \sum_{j=i}^{J} \pi_j \right)
\end{equation}
and hence
\begin{align*}
\label{dMtau}
\frac{\dd M}{\dd \tau_i} &= \sum_{j = 1}^{J} A_j^2 \frac{\dd \pi_j}{\dd \tau_i} - 2 \left( \sum_{j = 1}^{J} A_j \pi_j \right) \left( \sum_{j = 1}^{J} A_j \frac{\dd \pi_j}{\dd \tau_i} \right) \\
&= \alpha_i \Bigg( 
\left( \sum_{j=1}^{J} A_j^2 \pi_j \right) \left(\sum_{k = i}^J \pi_k \right) - \sum_{j=i}^{J} A_j^2 \pi_j \\
  &\qquad \quad - 2 \left( \sum_{j=1}^{J} A_j \pi_j \right) 
\left( 
\left(\sum_{j=1}^{J} A_j \pi_j \right) \left(\sum_{k = i}^J \pi_k \right) - \sum_{j=i}^{J} A_j \pi_j 
\right)
\Bigg). \numberthis
\end{align*}
Using Proposition \ref{tau_star_properties} (b), (\ref{dMtau}) evaluated in point $\tau^*$ is equal to 
\begin{equation}
\frac{\dd M}{\dd \tau_i} (\boldsymbol{\tau^*}, \boldsymbol{\alpha}) = \alpha_i A_J^2 \pi_J \pi_0 (\pi_J - \pi_0)
\end{equation}
Under suitable regularity conditions on $\Pi$, differentiation with respect to $\tau_i$ and integration with respect to $\theta$ may be interchanged so that we arrive at
\begin{align}
\frac{\dd \psi_0}{\dd \tau_i}(\boldsymbol{\tau^*}, \boldsymbol{\alpha}) &= \int \frac{\dd M}{\dd \tau_i}(\boldsymbol{\tau^*}, \boldsymbol{\alpha}) M(\boldsymbol{\tau^*}, \boldsymbol{\alpha})^{-1} \, \Pi \, \dd \theta = \alpha_i \int (\pi_J - \pi_0) \, \Pi \, \dd \theta \\
\frac{\dd \psi_{-1}}{\dd \tau_i}(\boldsymbol{\tau^*}, \boldsymbol{\alpha}) &= \int \frac{\dd M}{\dd \tau_i}(\boldsymbol{\tau^*}, \boldsymbol{\alpha}) M(\boldsymbol{\tau^*}, \boldsymbol{\alpha})^{-2} \, \Pi \, \dd \theta = \frac{\alpha_i}{A_J^2} \int \frac{(\pi_J - \pi_0)}{\pi_J \pi_0} \ \, \Pi \, \dd \theta
\end{align}
as
\begin{equation}
\label{Mtau_star}
M(\boldsymbol{\tau^*}, \boldsymbol{\alpha}) = A_k^2 \pi_J \pi_0
\end{equation} 
By symmetry of $\Pi$ and $\pi_j(. \, , \tau^*)$ (in the sense of Proposition \ref{tau_star_properties} (d)) we have
\begin{equation}
\int \pi_J^2 \pi_0 \, \Pi \, \dd \theta = \int \pi_J \pi_0^2 \, \Pi \, \dd \theta
\end{equation}
and hence $\frac{d \psi_0}{d \tau_i}(\boldsymbol{\tau^*}, \boldsymbol{\alpha}) = \frac{d \psi_{-1}}{d \tau_i}(\boldsymbol{\tau^*}, \boldsymbol{\alpha}) = 0$ for all $i \leq J$ due to symmetry. Since $\psi_v$ ($v \in \{-1, 0\}$) is concave \cite{firth1997}, $\boldsymbol{\tau^*}$ is a global maximum of the two design criteria within the set of one-point designs independent of the choice of $\boldsymbol{\alpha}$. 
\hfill $\Box$
\end{proof}


\begin{proof} \emph{(Lemma \ref{BOD_one_point_alpha})}
Although, we only need $J = 1$ for the proof, we compute the derivatives of $M$ with respect to $\alpha$ for arbitrary $J$. For $j \geq i $ the derivatives of $\pi_j$ after $\alpha_i$ equal
\begin{align*}
\frac{\dd \pi_j}{\dd \alpha_i} &= - (\theta - \tau_i) \left(\exp\left(\sum_{s=1}^j \alpha_s (\theta - \tau_s) \right) \left( \sum_{k = i}^J \exp\left(\sum_{s=1}^k \alpha_s (\theta - \tau_s) \right) - z \right) \right) / z^2 \\
&= -(\theta - \tau_i) \pi_j \left(\sum_{k = i}^J \pi_k - 1 \right) = 
-(\theta - \tau_i) \left( \pi_j \left(\sum_{k = i}^J \pi_k \right) - \pi_j \right) \\. \numberthis
\end{align*}
For $j < i$ the derivatives equal
\begin{align*}
\frac{\dd \pi_j}{\dd \alpha_i} &= -(\theta - \tau_i) \left(\exp\left(\sum_{s=1}^j \alpha_s (\theta - \tau_s) \right) \left(\sum_{k = i}^J \exp\left(\sum_{s=1}^k \alpha_s (\theta - \tau_s) \right) \right) \right) / z^2 \\
&= -(\theta - \tau_i) \pi_j \left( \sum_{k = i}^J \pi_k \right). \numberthis
\end{align*}

Thus, we find
\begin{align*}
\sum_{j = 1}^{J} \frac{\dd \pi_j}{\dd \alpha_i} &= -(\theta - \tau_i) \left( \left(\sum_{j=1}^{J} \pi_j \right) \left(\sum_{k = i}^J \pi_k \right) -\sum_{j=i}^{J} \pi_j \right)  \numberthis
\end{align*}
and hence 
\begin{align*}
\label{dMalpha}
\frac{\dd M}{\dd \alpha_i} &= \sum_{j = 1}^{J} \left(2 A_j\frac{\dd A_j}{\dd \alpha_i} \pi_j + A_j^2 \frac{\dd \pi_j}{\dd \alpha_i} \right) - 2 \left( \sum_{j = 1}^{J} A_j \pi_j \right) \left( \sum_{j = 1}^{J} \left(\frac{\dd A_j}{\dd \alpha_i} \pi_j + A_j \frac{\dd \pi_j}{\dd \alpha_i} \right) \right) \\
&= -(\theta - \tau_i) \Bigg( 
\left( \sum_{j=1}^{J} A_j^2 \pi_j \right) \left(\sum_{k = i}^J \pi_k \right) - \sum_{j=i}^{J} A_j^2 \pi_j \\
  &\qquad \quad - 2 \left( \sum_{j=1}^{J} A_j \pi_j \right) 
\left( 
\left(\sum_{j=1}^{J} A_j \pi_j \right) \left(\sum_{k = i}^J \pi_k \right) - \sum_{j=i}^{J} A_j \pi_j 
\right)
\Bigg) \\
 &\qquad \quad + \sum_{j = i}^{J} 2 A_j \pi_j -  
 2 \left( \sum_{j=1}^{J} A_j \pi_j \right) \left(\sum_{j=i}^{J} \pi_j \right).
 \numberthis
\end{align*}

Using Proposition \ref{tau_star_properties} (b), (\ref{dMalpha}) evaluated in point $\tau^*$ is equal to 
\begin{equation}
\frac{\dd M}{\dd \alpha_i} (\boldsymbol{\tau^*}, \boldsymbol{\alpha}) = A_J \pi_J \pi_0 \left(2 - A_J (\theta - \tau^*_i) (\pi_J - \pi_0) \right),
\end{equation}
which for $J = 1$ (dichotomous 2PL model; $A_J = \alpha \in \mathbb{R}_+$), can be written as
\begin{equation}
\frac{\dd M}{\dd \alpha} (0, \alpha) = \alpha \pi_1 \pi_0 \left(2 - \alpha \theta (\pi_1 - \pi_0) \right)
\end{equation}
Due to symmetry of $\pi$ in the sense of Proposition \ref{tau_star_properties} (d) and since $\Pi$ is symmetric, we have 
\begin{equation}
-\int \theta \pi_1 (\pi_1 \pi_0)^x \, \Pi \dd \theta = \int \theta \pi_0 (\pi_1 \pi_0)^x \, \Pi \dd \theta  
\end{equation}
for any $x \in \mathbb{Z}$. For symmetric $\Pi$, the integral $\int \pi_1 \pi_0 \, \Pi \, \dd \theta$ is monotonically decreasing in $\alpha$, whereas 
\begin{equation}
I_0(\alpha) := \int \theta \pi_1 \, \Pi \, \dd \theta
\end{equation}
is monotonically increasing, which both is immediately evident from the graph of $\pi_1$ and $\pi_0$ changing with $\alpha$. For $\psi_0$ we have
\begin{align}
\frac{\dd \psi_0}{\dd \alpha}(0, \alpha) &= \int \frac{\dd M}{\dd \alpha}(0, \alpha) M(0, \alpha)^{-1} \, \Pi \, \dd \theta = \frac{2}{\alpha} \int (1 - \alpha \theta \pi_1) \, \Pi \, \dd \theta
\end{align}
Since $I_0(\alpha)$ is monotonically increasing in $\alpha$, $\alpha I_0(\alpha)$ is monotonically increasing and unbounded. Furthermore, $\lim_{\alpha \rightarrow 0} \alpha I_0(\alpha) = 0$ and hence there exists a unique solution $\alpha^*_0$ of
\begin{equation}
\alpha \int \theta \pi_1 \, \Pi \, \dd \theta = 1,
\end{equation}
which (together with $\tau^*$) constitutes the optimal one-point design for $\psi_0$. For $\psi_{-1}$ we have
\begin{align}
\frac{\dd \psi_{-1}}{\dd \alpha}(0, \alpha) &= \int \frac{\dd M}{\dd \alpha}(0, \alpha) M(0, \alpha)^{-2} \, \Pi \, \dd \theta = \frac{2}{\alpha^3} \int \frac{(1 - \alpha \theta \pi_1)}{\pi_1 \pi_0} \, \Pi \, \dd \theta \\
\end{align}
The integral 
\begin{equation}
I_{-1}(\alpha) := \int \frac{(1 - \alpha \theta \pi_1)}{\pi_1 \pi_0} \, \Pi \, \dd \theta
\end{equation}
is monotonically decreasing in $\alpha$ for symmetric $\Pi$, which again can be inferred from the graph of $\pi_1$ and $\pi_0$ changing with $\alpha$. Since $\lim_{\alpha \rightarrow 0} I_{-1}(\alpha) = 4$ and $\lim_{\alpha \rightarrow \infty} I_{-1}(\alpha) = -\infty$, there exists a unique solution $\alpha^*_{-1}$ of
\begin{equation}
\alpha \int \frac{\theta}{\pi_0} \, \Pi \, \dd \theta = \int \frac{1}{\pi_1 \pi_0} \, \Pi \, \dd \theta,
\end{equation}
which (together with $\tau^*$) constitutes the optimal one-point design for $\psi_{-1}$. \hfill $\Box$
\end{proof}

\begin{proof} \emph{(Theorem \ref{BOD_PCM})}
For a symmetric weight distribution, the first derivatives of $\phi_0$ and $\phi_{-1}$ are zero in point $\tau^*$, which can be seen directly from the first derivatives of $\psi_0$ and $\psi_{-1}$ in $\tau^*$. For $i \geq n$ we compute the elements of the Hessian matrix of $M$ as
\begin{align*}
\label{d2Mtau}
\frac{\dd^2 M}{\dd \tau_i \tau_n} &= 
\alpha_i \alpha_n \Bigg( \left( \left( \sum_{j=1}^{J} A_j^2 \pi_j \right) \left(\sum_{k = n}^J \pi_k \right) - \left( \sum_{j=n}^{J} A_j^2 \pi_j \right) \right) \left(\sum_{k = i}^J \pi_k \right) \\
&\quad + \left(\sum_{j=1}^{J} A_j^2 \pi_j \right) \left( \left( \sum_{k = i}^J \pi_k \right) \left(\sum_{v = n}^J \pi_v \right) - \left( \sum_{k = i}^J \pi_k \right) \right) \\
&\quad - \left( \sum_{j = i}^{J} A_j^2 \pi_j \right) \left(\sum_{k = n}^J \pi_k \right) - \left( \sum_{j = i}^{J} A_j^2 \pi_j \right) \\
&\quad - 2 \left( \left( \sum_{j=1}^{J} A_j \pi_j \right) \left(\sum_{k = n}^J \pi_k \right) - \left( \sum_{j=n}^{J} A_j \pi_j \right) \right) 
\left( 
\left( \sum_{j=1}^{J} A_j \pi_j \right) \left(\sum_{k = i}^J \pi_k \right) - \left( \sum_{j=i}^{J} A_j \pi_j \right) 
\right) \\ 
&\quad - 2 \left( \sum_{j=1}^{J} A_j \pi_j \right) 
\Bigg(
\left( \left( \sum_{j=1}^{J} A_j \pi_j \right) \left(\sum_{k = n}^J \pi_k \right) - \left( \sum_{j=n}^{J} A_j \pi_j \right) \right)\left(\sum_{k = i}^J \pi_k \right) \\
&\qquad \qquad \qquad \qquad \quad + \left( \sum_{j=1}^{J} A_j \pi_j \right) \left( \left( \sum_{k = i}^J \pi_k \right) \left(\sum_{v = n}^J \pi_v \right) - \left( \sum_{k = i}^J \pi_k \right) \right) \\
&\qquad \qquad \qquad \qquad \quad - \left( \sum_{j = i}^{J} A_j \pi_j \right) \left(\sum_{k = n}^J \pi_k \right) - \left( \sum_{j = i}^{J} A_j \pi_j \right) \Bigg)
\Bigg). \numberthis
\end{align*}
This holds for $i < n$ as well due to symmetry of the Hessian matrix. Equation (\ref{d2Mtau}) can be slightly simplified, but remains too complicated to determine analytically under which conditions the Hessian matrix of $\phi_v$ ($v \in \{0, -1\}$) is positive definite for some Bayes optimal one-point design $\tau_{\Pi}$ for a weight distribution $\Pi$ of unspecified form. If, however, $\Pi$ is symmetric we know from Lemma \ref{BOD_one_point} that $\tau^*$ is the Bayes optimal one-point design and (\ref{d2Mtau}) remarkably simplifies to 
\begin{equation}
\frac{\dd^2 M}{\dd \tau_i \tau_n}(\boldsymbol{\tau^*}, \boldsymbol{\alpha}) = \alpha_i \alpha_n \pi_J \pi_0 A_J^2 (1-6 \pi_J \pi_0) = \alpha_i \alpha_n M(\boldsymbol{\tau^*}, \boldsymbol{\alpha})  (1-6 \pi_J \pi_0) 
\end{equation}
and hence
\begin{equation}
\label{dphi0tau}
\frac{\dd^2 \phi_0}{\dd \tau_i \tau_n}(\boldsymbol{\tau^*}, \boldsymbol{\alpha}) = \int M(\boldsymbol{\tau^*}, \boldsymbol{\alpha})^{-1} \frac{\dd^2 M}{\dd \tau_i \tau_n}(\boldsymbol{\tau^*}, \boldsymbol{\alpha}) \, \Pi \, \dd \theta = \alpha_i \alpha_n \int (1 - 6 \pi_J \pi_0) \, \Pi \, \dd \theta.
\end{equation}
Within the design space of the threshold vector $\tau$, this result is of limited use since the $\tau^*$ is at the border of the design region and the Hessian matrix is constant. However, considering that in $\tau^*$ the GPCM reduces to the dichotomous 2PL model, the necessary condition for $\tau^*$ to be Bayes optimal for $v = 0$ can be directly inferred from (\ref{dphi0tau}) as
\begin{equation}
\int (1 - 6 \pi_J \pi_0) \, \Pi \, \dd \theta \leq 0
\end{equation}
or equivalently
\begin{equation}
\int \pi_J \pi_0 \, \Pi \, \dd \theta \geq \frac{1}{6}.
\end{equation}
Similarly, when denoting $C := (\int M(\boldsymbol{\tau^*}, \boldsymbol{\alpha})^{-1} \, \Pi \, \dd \theta)^{-1}$, we obtain
\begin{align*}
\label{dphim1tau}
\frac{\dd^2 \phi_{-1}}{\dd \tau_i \tau_n}(\boldsymbol{\tau^*}, \boldsymbol{\alpha}) &= 
C \int M(\boldsymbol{\tau^*}, \boldsymbol{\alpha})^{-2} \frac{\dd^2 M}{d \tau_i \tau_n}(\boldsymbol{\tau^*}, \boldsymbol{\alpha}) \, \Pi \, \dd \theta \\
&= \alpha_i \alpha_n C \int (\pi_J \pi_0 J^2)^{-1} (1 - 6 \pi_J \pi_0) \, \Pi \, \dd \theta \\
&= \alpha_i \alpha_n \left( C \int M(\boldsymbol{\tau^*}, \boldsymbol{\alpha})^{-1} \, \Pi \, \dd \theta - C \frac{6}{J^2} \int \, \Pi \, \dd \theta \right) \\
&= \alpha_i \alpha_n \left( 1 - C \frac{6}{J^2} \right).
\end{align*}
Thus, the necessary condition for $\tau^*$ to be Bayes optimal for $v = -1$ becomes
\begin{equation}
6 \geq \int (\pi_J \pi_0)^{-1} \, \Pi \, \dd \theta.
\end{equation}
\hfill $\Box$
\end{proof}

\begin{proof} \emph{(Theorem \ref{BOD_PCM_symmetric})}
This can be proved in the same way as Theorem 2 in \citeA{grasshoff2012}, so we do not spell out the details here. \hfill $\Box$
\end{proof}


\vspace{\fill}\pagebreak





\bibliography{citations}{}
\bibliographystyle{apacite}


%




\end{document}